© 2021 г.     **А.А. БОБЦОВ,** д-р техн. наук (bobtsov@mail.ru)
**Н.А. НИКОЛАЕВ,** к-т техн. наук (nikona@yandex.ru)
**О.В. ОСЬКИНА,** аспирант (olga.oskina1996@gmail.com)
**С. И. НИЗОВЦЕВ,** аспирант (nizovtsev.si@gmail.com)
**(Университет ИТМО г. Санкт-Петербург)**

# ИДЕНТИФИКАЦИЯ НЕСТАЦИОНАРНОЙ ЧАСТОТЫ НЕ ЗАШУМЛЕННОГО СИНУСОИДАЛЬНОГО СИГНАЛА

*Аннотация.* Рассматривается новый алгоритм оценивания нестационарной частоты незашумленного синусоидального сигнала. Предполагается, что амплитуда и частота синусоидального сигнала неизвестные функции времени, но являются решениями линейных стационарных дифференциальных уравнений с известными параметрами. Поставленная задача решается с использованием градиентных алгоритмов настройки на базе линейного регрессионного уравнения, полученного путем параметризации исходного нелинейного по параметрам синусоидального сигнала. Представленный в статье пример и результаты компьютерного моделирования иллюстрируют работоспособность предлагаемого алгоритма, а также поясняют процедуру его синтеза.

*Ключевые слова*: идентификация параметров, нестационарные системы, синусоидальные сигналы

## 1. Введение.

В статье рассматривается задача идентификации неизвестной нестационарной частоты незашумленного синусоидального сигнала вида

(1)     $y(t) = \alpha(t) \sin(\omega(t) + \varphi)$,

где $y(t)$ – сигнал, доступный прямому измерению, $\alpha(t)$ – неизвестная нестационарная амплитуда, $\omega(t)$ – неизвестный нестационарный параметр, $\varphi$ – неизвестный постоянный фазовый сдвиг.

В случае стационарных параметров $\alpha$ и $\omega$ задача оценивания неизвестной частоты по измерению сигнала $y(t)$ хорошо изучена, и ей посвящено большое число публикаций (см., например, [1 – 9]). Однако такое допущение может не выполняться при решении реальных инженерных задач. Если источником возмущения является работа электрического привода, то частота пропорциональна скорости вращения электрической машины, а следовательно, изменяется при разгоне или торможении. Для парирования данного эффекта были предложены методы, опубликованные в статьях [10 – 12], в которых частота синусоидального сигнала описывается полиномом времени произвольного порядка. Тем не менее в работах [10 – 12] допускалось, что амплитуда $\alpha$ является неизвестным числом. В работах [13 – 15] ограничения на стационарность амплитуды были сняты, но частота синусоидального сигнала является постоянной. Таким образом, на сколько известно авторам, тематика идентификации нестационарной частоты при условии нестационарной амплитуды является открытой, что в свою очередь, мотивирует проведение новых исследований.

## 2. Постановка задачи

По измерению сигнала (1) ставится задача синтеза алгоритма идентификации параметра $\omega$ при следующих допущениях.

*Допущение 1.* Функция $\omega(t)$ является выходом линейного генератора

(2) $\qquad \omega = h^T \xi,$

(3) $\qquad \dot{\xi} = \Gamma \xi,$

где параметры вектора $h$ и матрицы $\Gamma$ являются известными числами, но начальные условия вектора $\xi$ неизвестны.

*Допущение 2.* Функция $\alpha(t)$ является выходом линейного генератора

(4) $\quad \alpha = r^T \eta,$

(5) $\quad \dot{\eta} = G\eta,$

где параметры вектора $r$ и матрицы $G$ являются известными числами, но начальные условия вектора $\eta$ неизвестны.

## 3. Основной результат

Для синтеза алгоритма оценивания неизвестной нестационарной частоты синусоидального сигнала (1) будем использовать обобщенный подход к синтезу наблюдателей, основанный на оценке параметров (GPEBO – Generalized parameter estimation-based observers) [16]. Первым шагом является преобразование исходной нелинейной модели вида (1) к линейной регрессии, параметры которой, в свою очередь, могут быть оценены с использованием, например, стандартного градиентного подхода (см., например [17]) или подхода динамического расширения регрессора и смешивания (DREM – Dynamic regressor extension and mixing) [18].

Для преобразования исходной нелинейной по параметру $\omega$ модели (1) продифференцируем сигнал $y(t) = \alpha(t)\sin(\omega(t) + \varphi)$ два раза. Для первой производной (1) имеем

$$\dot{y} = \dot{\alpha}\sin(\omega + \varphi) + \alpha\cos(\omega + \varphi)\dot{\omega}$$

или в более удобном виде

(6) $\quad \dot{y} = \frac{\dot{\alpha}}{\alpha}y + \alpha\dot{\omega}\cos(\omega + \varphi).$

Для второй производной получаем

(7) $\quad \ddot{y} = \ddot{\alpha}\sin(\omega + \varphi) + \dot{\alpha}\cos(\omega + \varphi)\dot{\omega} + \dot{\alpha}\cos(\omega + \varphi)\dot{\omega} -$

$-\alpha\sin(\omega + \varphi)\dot{\omega}^2 + \alpha\cos(\omega + \varphi)\ddot{\omega} =$

$= \ddot{\alpha}\sin(\omega + \varphi) + 2\dot{\alpha}\cos(\omega + \varphi)\dot{\omega} - \underbrace{\alpha\sin(\omega + \varphi)}_{y}\dot{\omega}^2 + \alpha\cos(\omega + \varphi)\ddot{\omega} =$

$= \ddot{\alpha}\sin(\omega + \varphi) + 2\dot{\alpha}\dot{\omega}\cos(\omega + \varphi) - \dot{\omega}^2 y + \alpha\ddot{\omega}\cos(\omega + \varphi).$

Подставим в (7) выражение для $\cos(\omega + \varphi) = \frac{\alpha\dot{y} - \dot{\alpha}y}{\alpha^2\dot{\omega}}$, которое можно получить из (6)

(8) $\quad \ddot{y} = \ddot{\alpha}\sin(\omega + \varphi) + 2\dot{\alpha}\dot{\omega}\frac{(\alpha\dot{y} - \dot{\alpha}y)}{\alpha^2\dot{\omega}} - \dot{\omega}^2 y + \alpha\ddot{\omega}\frac{(\alpha\dot{y} - \dot{\alpha}y)}{\alpha^2\dot{\omega}} =$

$= \ddot{\alpha}\underbrace{\sin(\omega + \varphi)}_{\frac{y}{\alpha}} + 2\dot{\alpha}\frac{(\alpha\dot{y} - \dot{\alpha}y)}{\alpha^2} - \dot{\omega}^2 y + \ddot{\omega}\frac{(\alpha\dot{y} - \dot{\alpha}y)}{\alpha\dot{\omega}} =$

$= \frac{\ddot{\alpha}}{\alpha}y + \frac{2\dot{\alpha}\alpha}{\alpha^2}\dot{y} - \frac{2\dot{\alpha}^2}{\alpha^2}y - \dot{\omega}^2 y + \frac{\ddot{\omega}\alpha}{\alpha\dot{\omega}}\dot{y} - \frac{\ddot{\omega}\dot{\alpha}}{\alpha\dot{\omega}}y =$

$= \left(\frac{\ddot{\alpha}}{\alpha} - \frac{2\dot{\alpha}^2}{\alpha^2} - \dot{\omega}^2 - \frac{\dot{\alpha}\ddot{\omega}}{\alpha\dot{\omega}}\right)y + \left(\frac{2\dot{\alpha}}{\alpha} + \frac{\ddot{\omega}}{\dot{\omega}}\right)\dot{y}.$

После выполнения простейших математических преобразований имеем

(9) $\quad \alpha^2\dot{\omega}\ddot{y} = (\alpha\ddot{\alpha}\dot{\omega} - 2\dot{\alpha}^2\dot{\omega} - \alpha^2\dot{\omega}^3 - \alpha\dot{\alpha}\ddot{\omega})y + (2\alpha\dot{\alpha}\dot{\omega} + \alpha^2\ddot{\omega})\dot{y}.$

Для вывода линейной регрессионной модели дважды профильтруем левую и правую части выражения (9) с использованием апериодического звена первого порядка $\frac{1}{(p+1)}$. При выполнении математических преобразований будем использовать лемму о перестановке ("swapping lemma" [17]).

Шаг 1. Пропустим левую часть (9) через фильтр $\frac{1}{(p+1)}$

$$\frac{1}{p+1}[(\alpha^2\dot{\omega})\ddot{y}] = (\alpha^2\dot{\omega})\frac{1}{p+1}[\ddot{y}] - \frac{1}{p+1}\left[p(\alpha^2\dot{\omega})\frac{1}{p+1}[\ddot{y}]\right] =$$

$$= (\alpha^2\dot{\omega})\underbrace{\frac{p}{p+1}[\dot{y}]}_{\dot{q}_1} - \frac{1}{p+1}\left[(2\alpha\dot{\alpha}\dot{\omega} + \alpha^2\ddot{\omega})\underbrace{\frac{p}{p+1}[\dot{y}]}_{\dot{q}_1}\right] =$$

$$= \alpha^2\dot{\omega}\dot{q}_1 - \frac{1}{p+1}[(2\alpha\dot{\alpha}\dot{\omega} + \alpha^2\ddot{\omega})\dot{q}_1] =$$

$$= \alpha^2\dot{\omega}\dot{q}_1 - \left((2\alpha\dot{\alpha}\dot{\omega} + \alpha^2\ddot{\omega})\frac{1}{p+1}[\dot{q}_1] - \frac{1}{p+1}\left[p(2\alpha\dot{\alpha}\dot{\omega} + \alpha^2\ddot{\omega})\frac{1}{p+1}[\dot{q}_1]\right]\right) =$$

$$= \alpha^2\dot{\omega}\dot{q}_1 - \left((2\alpha\dot{\alpha}\dot{\omega} + \alpha^2\ddot{\omega})\underbrace{\frac{p}{p+1}[q_1]}_{q_2} - \frac{1}{p+1}\left[p(2\alpha\dot{\alpha}\dot{\omega} + \alpha^2\ddot{\omega})\underbrace{\frac{p}{p+1}[q_1]}_{q_2}\right]\right) =$$

$$= \alpha^2\dot{\omega}\dot{q}_1 - \left((2\alpha\dot{\alpha}\dot{\omega} + \alpha^2\ddot{\omega})q_2 - \frac{1}{p+1}\left[\left(p(2\alpha\dot{\alpha}\dot{\omega}) + p(\alpha^2\ddot{\omega})\right)q_2\right]\right) =$$

$$= \alpha^2\dot{\omega}\dot{q}_1 - (2\alpha\dot{\alpha}\dot{\omega} + \alpha^2\ddot{\omega})q_2 + \frac{1}{p+1}[(2\dot{\alpha}^2\dot{\omega} + 2\alpha\ddot{\alpha}\dot{\omega} + 4\alpha\dot{\alpha}\ddot{\omega} + \alpha^2\dddot{\omega})q_2],$$

где $p = d/dt$ – оператор дифференцирования, а функции $q_1$ и $q_2$ имеют вид

(10) $\begin{cases} q_1 = \frac{p}{p+1}y \\ q_2 = \frac{p}{p+1}q_1 \end{cases}$,

таким образом имеем

(11) $\frac{1}{p+1}[(\alpha^2\dot{\omega})\ddot{y}] = \alpha^2\dot{\omega}\dot{q}_1 - (2\alpha\dot{\alpha}\dot{\omega} + \alpha^2\ddot{\omega})q_2 +$

$+ \frac{1}{p+1}[(2\dot{\alpha}^2\dot{\omega} + 2\alpha\ddot{\alpha}\dot{\omega} + 4\alpha\dot{\alpha}\ddot{\omega} + \alpha^2\dddot{\omega})q_2].$

Шаг 2. Пропустим (11) через фильтр $\frac{1}{(p+1)}$

$$\frac{1}{(p+1)^2}[(\alpha^2\dot{\omega})\ddot{y}] = \frac{1}{(p+1)}[(\alpha^2\dot{\omega})\dot{q}_1] - \frac{1}{(p+1)}[(2\alpha\dot{\alpha}\dot{\omega} + \alpha^2\ddot{\omega})q_2] +$$

$$+ \frac{1}{(p+1)^2}[(2\dot{\alpha}^2\dot{\omega} + 2\alpha\ddot{\alpha}\dot{\omega} + 4\alpha\dot{\alpha}\ddot{\omega} + \alpha^2\dddot{\omega})q_2] =$$

$$= (\alpha^2\dot{\omega})\underbrace{\frac{p}{(p+1)}[q_1]}_{q_2} - \frac{1}{(p+1)}[(2\alpha\dot{\alpha}\dot{\omega} + \alpha^2\ddot{\omega})q_2] - \frac{1}{(p+1)}[(2\alpha\dot{\alpha}\dot{\omega} + \alpha^2\ddot{\omega})q_2] +$$

$$+ \frac{1}{(p+1)^2}[(2\dot{\alpha}^2\dot{\omega} + 2\alpha\ddot{\alpha}\dot{\omega} + 4\alpha\dot{\alpha}\ddot{\omega} + \alpha^2\dddot{\omega})q_2] =$$

$$= \alpha^2\dot{\omega}q_2 - \frac{2}{(p+1)}[(2\alpha\dot{\alpha}\dot{\omega} + \alpha^2\ddot{\omega})q_2] +$$

$$+ \frac{1}{(p+1)^2}[(2\dot{\alpha}^2\dot{\omega} + 2\alpha\ddot{\alpha}\dot{\omega} + 4\alpha\dot{\alpha}\ddot{\omega} + \alpha^2\dddot{\omega})q_2]$$

таким обратом, для дважды профильтрованной левой части (9) имеем

(12) $\frac{1}{(p+1)^2}[(\alpha^2\dot{\omega})\ddot{y}] = \alpha^2\dot{\omega}q_2 - \frac{2}{(p+1)}[(2\alpha\dot{\alpha}\dot{\omega} + \alpha^2\ddot{\omega})q_2] +$

$+ \frac{1}{(p+1)^2}[(2\dot{\alpha}^2\dot{\omega} + 2\alpha\ddot{\alpha}\dot{\omega} + 4\alpha\dot{\alpha}\ddot{\omega} + \alpha^2\dddot{\omega})q_2].$

Шаг 3. Пропустим дважды правую часть (9) через фильтр $\frac{1}{(p+1)}$

$$\frac{1}{(p+1)^2}[(\alpha\ddot{\alpha}\dot{\omega} - 2\dot{\alpha}^2\dot{\omega} - \alpha^2\dot{\omega}^3 - \alpha\dot{\alpha}\ddot{\omega})y + (2\alpha\dot{\alpha}\dot{\omega} + \alpha^2\ddot{\omega})\dot{y}] =$$

$$= \frac{1}{(p+1)}\left[\frac{1}{(p+1)}[(\alpha\ddot{\alpha}\dot{\omega} - 2\dot{\alpha}^2\dot{\omega} - \alpha^2\dot{\omega}^3 - \alpha\dot{\alpha}\ddot{\omega})y + (2\alpha\dot{\alpha}\dot{\omega} + \alpha^2\ddot{\omega})\dot{y}]\right] =$$

$$= \frac{1}{(p+1)^2}[(\alpha\ddot{\alpha}\dot{\omega} - 2\dot{\alpha}^2\dot{\omega} - \alpha^2\dot{\omega}^3 - \alpha\dot{\alpha}\ddot{\omega})y] +$$

$$+ \frac{1}{(p+1)}\left[(2\alpha\dot{\alpha}\dot{\omega} + \alpha^2\ddot{\omega})\frac{1}{(p+1)}[\dot{y}] - \frac{1}{(p+1)}\left[p(2\alpha\dot{\alpha}\dot{\omega} + \alpha^2\ddot{\omega})\frac{1}{p+1}[\dot{y}]\right]\right] =$$

$$= \frac{1}{(p+1)^2}[(\alpha\ddot{\alpha}\dot{\omega} - 2\dot{\alpha}^2\dot{\omega} - \alpha^2\dot{\omega}^3 - \alpha\dot{\alpha}\ddot{\omega})y] +$$

$$+ \frac{1}{(p+1)}\left[(2\alpha\dot{\alpha}\dot{\omega} + \alpha^2\ddot{\omega})\underbrace{\frac{p}{(p+1)}[y]}_{q_1} - \frac{1}{(p+1)}\left[p(2\alpha\dot{\alpha}\dot{\omega} + \alpha^2\ddot{\omega})\underbrace{\frac{p}{(p+1)}[y]}_{q_1}\right]\right] =$$

$$= \frac{1}{(p+1)^2}[(\alpha\ddot{\alpha}\dot{\omega} - 2\dot{\alpha}^2\dot{\omega} - \alpha^2\dot{\omega}^3 - \alpha\dot{\alpha}\ddot{\omega})y] +$$

$$+ \frac{1}{(p+1)}\left[(2\alpha\dot{\alpha}\dot{\omega} + \alpha^2\ddot{\omega})q_1 - \frac{1}{(p+1)}[p(2\alpha\dot{\alpha}\dot{\omega} + \alpha^2\ddot{\omega})\ q_1]\right]$$

$$= \frac{1}{(p+1)^2}[(\alpha\ddot{\alpha}\dot{\omega} - 2\dot{\alpha}^2\dot{\omega} - \alpha^2\dot{\omega}^3 - \alpha\dot{\alpha}\ddot{\omega})y] +$$

$$+ \frac{1}{(p+1)}\left[(2\alpha\dot{\alpha}\dot{\omega} + \alpha^2\ddot{\omega})q_1 - \frac{1}{(p+1)}[(2\dot{\alpha}^2\dot{\omega} + 2\alpha\ddot{\alpha}\dot{\omega} + 4\alpha\dot{\alpha}\ddot{\omega} + \alpha^2\dddot{\omega})q_1]\right] =$$

$$= \frac{1}{(p+1)^2}[(\alpha\ddot{\alpha}\dot{\omega} - 2\dot{\alpha}^2\dot{\omega} - \alpha^2\dot{\omega}^3 - \alpha\dot{\alpha}\ddot{\omega})y] +$$

$$+ \frac{1}{(p+1)}[(2\alpha\dot{\alpha}\dot{\omega} + \alpha^2\ddot{\omega})q_1] - \frac{1}{(p+1)^2}[(2\dot{\alpha}^2\dot{\omega} + 2\alpha\ddot{\alpha}\dot{\omega} + 4\alpha\dot{\alpha}\ddot{\omega} + \alpha^2\dddot{\omega})q_1]$$

таким образом, после выполнения преобразований, имеем

(13) $\quad \frac{1}{(p+1)^2}[(\alpha\ddot{\alpha}\dot{\omega} - 2\dot{\alpha}^2\dot{\omega} - \alpha^2\dot{\omega}^3 - \alpha\dot{\alpha}\ddot{\omega})y + (2\alpha\dot{\alpha}\dot{\omega} + \alpha^2\ddot{\omega})\dot{y}] =$

$$= \frac{1}{(p+1)^2}[(\alpha\ddot{\alpha}\dot{\omega} - 2\dot{\alpha}^2\dot{\omega} - \alpha^2\dot{\omega}^3 - \alpha\dot{\alpha}\ddot{\omega})y] +$$

$$+ \frac{1}{(p+1)}[(2\alpha\dot{\alpha}\dot{\omega} + \alpha^2\ddot{\omega})q_1] - \frac{1}{(p+1)^2}[(2\dot{\alpha}^2\dot{\omega} + 2\alpha\ddot{\alpha}\dot{\omega} + 4\alpha\dot{\alpha}\ddot{\omega} + \alpha^2\dddot{\omega})q_1]$$

После объединения (12) и (13), для дважды профильтрованного уравнения (9) имеем

(14) $\quad \alpha^2\dot{\omega}q_2 - \frac{2}{(p+1)}[(2\alpha\dot{\alpha}\dot{\omega} + \alpha^2\ddot{\omega})q_2] + \frac{1}{(p+1)^2}[(2\dot{\alpha}^2\dot{\omega} + 2\alpha\ddot{\alpha}\dot{\omega} + 4\alpha\dot{\alpha}\ddot{\omega} + \alpha^2\dddot{\omega})q_2] =$

$$= \frac{1}{(p+1)^2}[(\alpha\ddot{\alpha}\dot{\omega} - 2\dot{\alpha}^2\dot{\omega} - \alpha^2\dot{\omega}^3 - \alpha\dot{\alpha}\ddot{\omega})y] +$$

$$+ \frac{1}{(p+1)}[(2\alpha\dot{\alpha}\dot{\omega} + \alpha^2\ddot{\omega})q_1] - \frac{1}{(p+1)^2}[(2\dot{\alpha}^2\dot{\omega} + 2\alpha\ddot{\alpha}\dot{\omega} + 4\alpha\dot{\alpha}\ddot{\omega} + \alpha^2\dddot{\omega})q_1].$$

Хорошо известно, что в силу допущений 1 и 2 функции $\alpha(t)$ и $\omega(t)$ могут быть записаны следующим образом

(15) $\quad \omega = h^T e^{\Gamma t}\, \xi_0,$

(16) $\quad \alpha = r^T e^{G t}\, \eta_0,$

где $\xi_0$ и $\eta_0$ – векторы неизвестных постоянных параметров.

Из (15) и (16) легко получить

$\dot{\omega} = h^T \Gamma e^{\Gamma t}\, \xi_0,\; \ddot{\omega} = h^T \Gamma^2 e^{\Gamma t}\, \xi_0,\; \dddot{\omega} = h^T \Gamma^3 e^{\Gamma t}\, \xi_0,\; \dot{\alpha} = r^T G e^{G t}\eta_0$ и $\ddot{\alpha} = r^T G^2 e^{G t}\eta_0.$

Откуда для уравнения (14) получаем линейную регрессионную модель относительно векторов неизвестных постоянных параметров $\xi_0$ и $\eta_0$.

## 4. Пример.

Для более наглядного понимания процедуры синтеза алгоритма оценивания функции $\omega(t)$ рассмотрим пример. Прежде всего предположим, что $\alpha$ – неизвестный постоянный параметр. Также будем допускать, что частота $\omega(t)$ изменяется по гармоническому закону, то есть формируется с помощью автономного генератора вида

(17) $\quad \omega = h^T \xi,$

(18) $\quad \dot{\xi} = \Gamma \xi,$

где $\Gamma = \begin{bmatrix} 0 & 1 \\ -\gamma & 0 \end{bmatrix}, h^T = [1 \quad 0].$

В соответствии с принятым допущением относительно неизвестных параметров, выражение (14) можно упростить и записать в виде

(19) $\quad \dot{\omega} q_2 = \frac{1}{(p+1)}[\ddot{\omega}(2q_2 + q_1)] - \frac{1}{(p+1)^2}[\ddot{\omega}(q_2 + q_1) + \dot{\omega}^3 y].$

Введем вспомогательную систему вида (см. [16])

(20) $\quad \hat{\omega} = h^T \xi_\omega,$

(21) $\quad \dot{\xi}_\omega = \Gamma \xi_\omega.$

Рассмотрим уравнение ошибки вида

(22) $\quad \varepsilon = \xi_\omega - \xi,$

тогда для производной от (22) имеем

(23) $\quad \dot{\varepsilon} = \Gamma \varepsilon.$

Решение дифференциального уравнения (18) имеет вид (см., например, [19])

$\varepsilon = e^{\Gamma t}\varepsilon(0) = \Phi \varepsilon(0) = \Phi \Theta,$

где $\Phi$ – фундаментальная матрица, $\dot{\Phi} = \Gamma\Phi$, $\Phi = I_{2\times 2}$, $\Theta = \begin{bmatrix} \theta_1 \\ \theta_2 \end{bmatrix}$ – неизвестные постоянные параметры, которые необходимо найти.

В случае, если начальные условия системы (20), (21) нулевые, то неизвестный вектор $\Theta$ является вектором начальных условий системы (17), (18)

$$\varepsilon(0) = \xi_\omega(0) - \xi_0 = -\xi_0.$$

Найдем производные сигнала (20)

(24) $\quad \widehat{\omega} = h^T \Phi \Theta,$

(25) $\quad \dot{\widehat{\omega}} = h^T \dot{\Phi} \Theta = h^T \Gamma \Phi \Theta = \Phi_{21} \theta_1 + \Phi_{22} \theta_2,$

(26) $\quad \ddot{\widehat{\omega}} = h^T \Gamma^2 \Phi \Theta = -\gamma h^T \Phi \Theta = \gamma(\Phi_{11} \theta_1 + \Phi_{12} \theta_2),$

(27) $\quad \dddot{\widehat{\omega}} = h^T \Gamma^3 \Phi \Theta = -\gamma h^T \Gamma \Phi \Theta = \gamma(\Phi_{21} \theta_1 + \Phi_{22} \theta_2).$

Подставляя выражения (24) – (26) в (19), получаем

(28) $\quad (\Phi_{21}\theta_1 + \Phi_{22}\theta_2)q_2 = \frac{-\gamma}{(p+1)}[(\Phi_{11}\theta_1 + \Phi_{12}\theta_2)(2q_2 + q_1)] +$

$\qquad + \frac{1}{(p+1)^2}[\gamma(\Phi_{21}\theta_1 + \Phi_{22}\theta_2)(q_2 + q_1) - (\Phi_{21}\theta_1 + \Phi_{22}\theta_2)^3 y].$

Таким образом регрессионная модель принимает вид

(29) $\quad z = m^T \Theta,$

где $z = \Phi_{21}q_2 + \frac{\gamma}{(p+1)}\Phi_{11}(2q_2 + q_1) - \frac{1}{(p+1)^2}\gamma\Phi_{21}(q_2 + q_1),$

$$m = \begin{bmatrix} -\Phi_{22}q_2 - \frac{1}{(p+1)}\gamma\Phi_{12}(2q_2 + q_1) + \frac{1}{(p+1)^2}\gamma\Phi_{22}(q_2 + q_1) \\ -\frac{1}{(p+1)^2}\Phi_{21}^3 y \\ -\frac{1}{(p+1)^2}3\Phi_{21}^2\Phi_{22} y \\ -\frac{1}{(p+1)^2}3\Phi_{21}\Phi_{22}^2 y \\ -\frac{1}{(p+1)^2}\Phi_{22}^3 y \end{bmatrix}, \Theta = \begin{bmatrix} \frac{\theta_2}{\theta_1} \\ \theta_1^2 \\ \theta_1\theta_2 \\ \theta_2^2 \\ \frac{\theta_2^3}{\theta_1} \end{bmatrix}.$$

Для идентификации постоянных неизвестных параметром модели (29) применим метод DREM. Для этого умножим левую и правую часть (29) слева на $m$

(30) $\quad mz = mm^T \Theta,$

Применим к (30) вспомогательный линейный фильтр $\frac{\lambda}{p+\lambda}$, тогда регрессионная модель (30) в новом базисе примет вид

(31) $\quad Y = \Omega\Theta,$

где $Y$ и $\Omega$ являются решением дифференциальных уравнений

(32) $\quad \dot{Y} = -\lambda Y + \lambda mz,$

(33) $\quad \dot{\Omega} = -\lambda\Omega + \lambda mm^T.$

Преобразуем регрессионную модель (31) к следующему виду

(34) $\quad \Upsilon = \Delta \Theta,$

где $\Upsilon = adj\Omega Y$, $adj\Omega$ – присоединенная матрица для $\Omega$, $\Delta = det\Omega$.

Оценку параметров регрессионной модели (34) выполним с помощью стандартного градиентного алгоритма идентификации вида

$$\dot{\hat{\Theta}} = -\beta\Delta(\Delta\hat{\Theta} - \Upsilon).$$

Проведем компьютерное моделирование для разных начальных значений $\xi_0$ и $\gamma$. На рисунке 1 приведены результаты моделирования для случая $\theta_1 = 2$, $\theta_2 = 1$, $\lambda = 1$, $\gamma = 4$ и $\beta = 10^{23}$. На рисунке 2 приведены результаты моделирования $\omega(t)$ и $\hat{\omega}(t)$ при $\theta_1 = 2$, $\theta_2 = 1$, $\lambda = 1$, $\gamma = 4$ и $\beta = 10^{23}$. На рисунке 3 приведены результаты моделирования $\omega(t)$ и $\hat{\omega}(t)$ при $\theta_1 = 4$, $\theta_2 = 2$, $\lambda = 1$, $\gamma = 4$ и $\beta = 10^{23}$. На рисунке 4 приведены графики ошибок $\hat{\omega}(t) - \omega(t)$ при $\lambda = 1$, $\gamma = 4$ и $\beta = 10^{23}$ для случаев $\theta_1 = 2$, $\theta_2 = 1$ и $\theta_1 = 4$, $\theta_2 = 2$. На рисунке 4 приведены графики ошибок $\hat{\omega}(t) - \omega(t)$ при $\theta_1 = 4$, $\theta_2 = 2$, $\lambda = 1$, $\gamma = 4$ и $\beta = 10^{23}$ для случаев $\gamma = 1$ и $\gamma = 4$. На рисунках 6 и 7 приведен график доступного прямым измерениям сигнала $y(t)$ при $\theta_1 = 2$, $\theta_2 = 1$, для $\gamma = 1$ и $\gamma = 4$ соответственно.

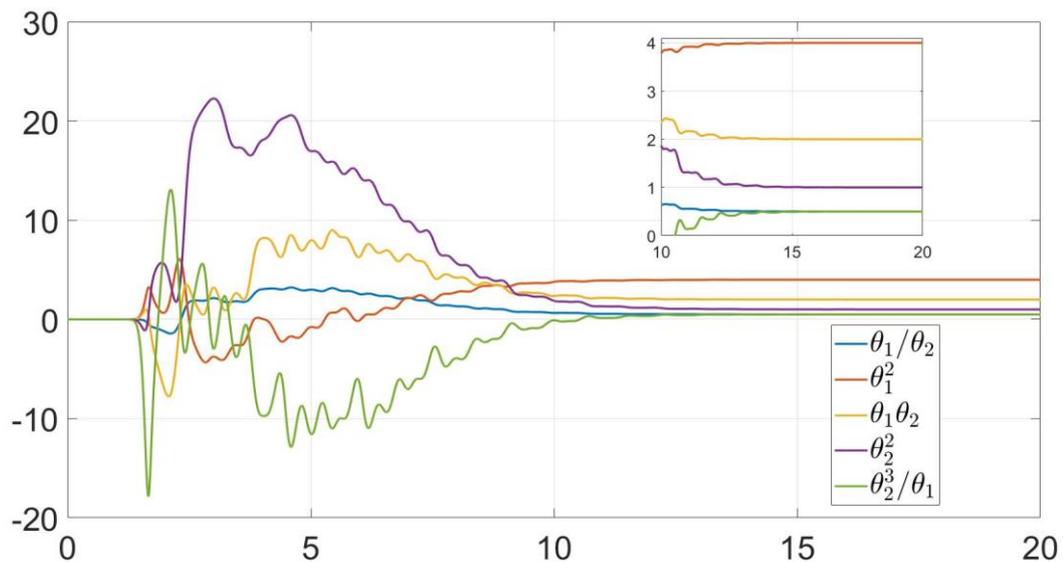

Рис. 1 – Результаты моделирования оценки неизвестных параметров регрессионной модели (29)

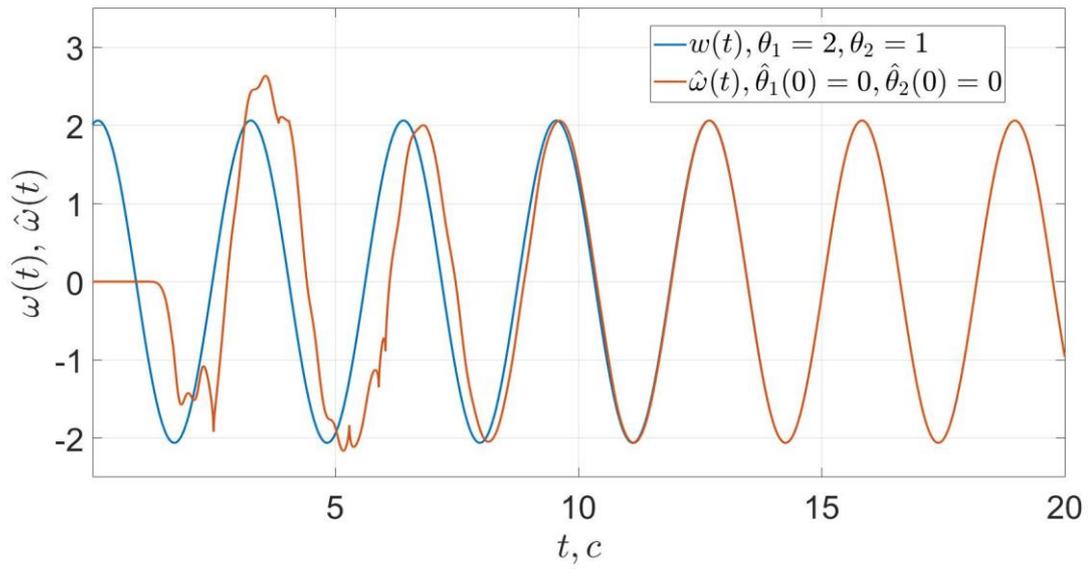

Рис. 2 – Графики $\omega(t)$ и $\hat{\omega}(t)$ для начальных условий $\theta_1 = 2, \theta_2 = 1$

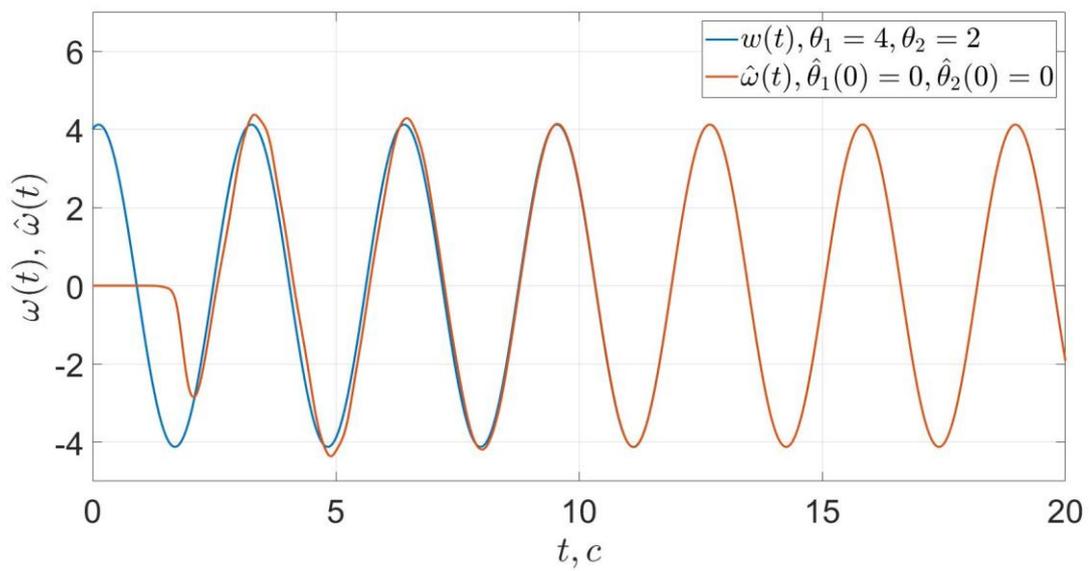

Рис. 3 – Графики $\omega(t)$ и $\hat{\omega}(t)$ для начальных условий $\theta_1 = 4, \theta_2 = 2$

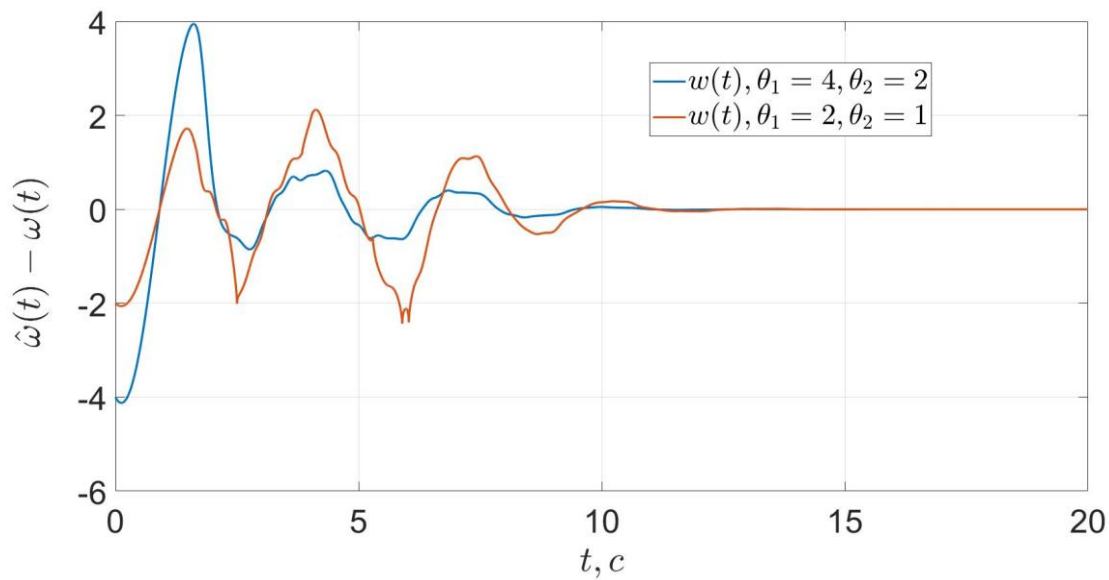

Рис. 4 – Графики ошибок при различных начальных условиях $\omega(t)$

Из результатов компьютерного моделирования следует сходимость оценки $\widehat{\omega}(t)$ к $\omega(t)$ для различных значений $\xi_0$ и $\gamma$.

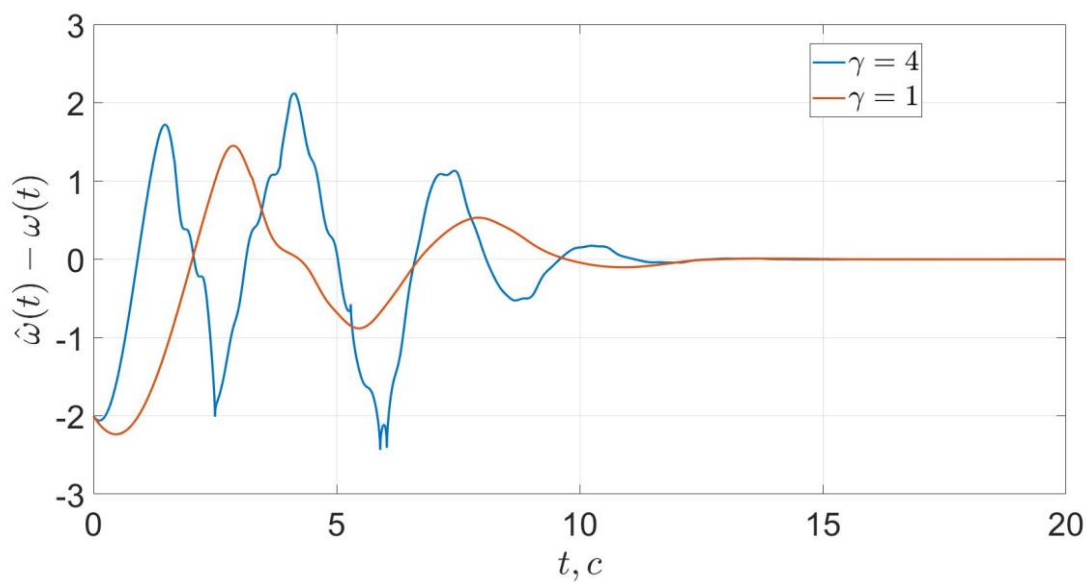

Рис. 5 – Графики ошибки для различных частот параметра $\omega(t)$

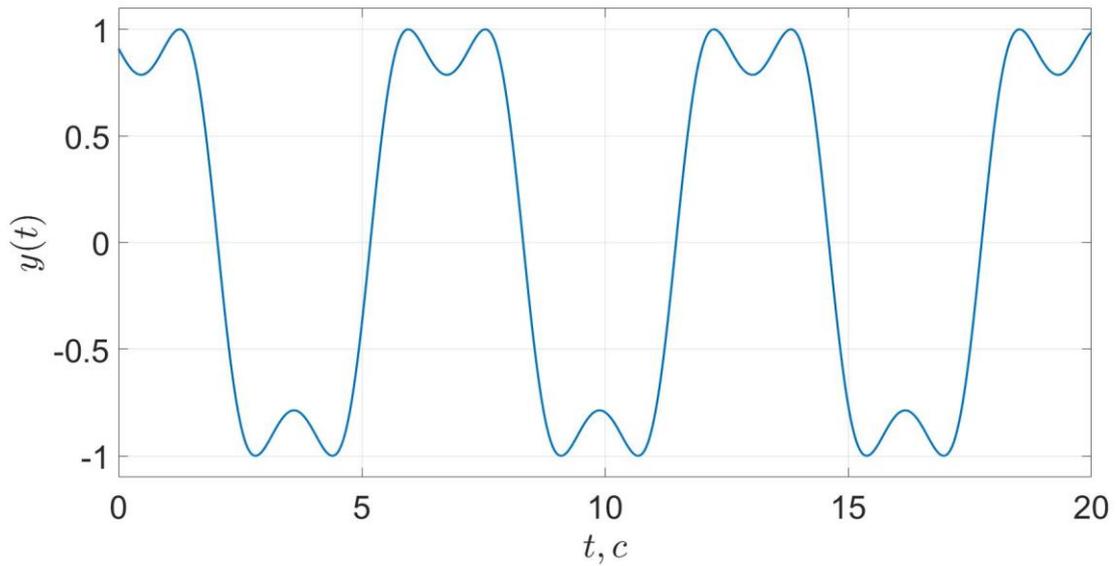

Рис. 6 – График сигнала $y(t)$ при $\theta_1 = 2, \theta_2 = 1, \gamma = 1$.

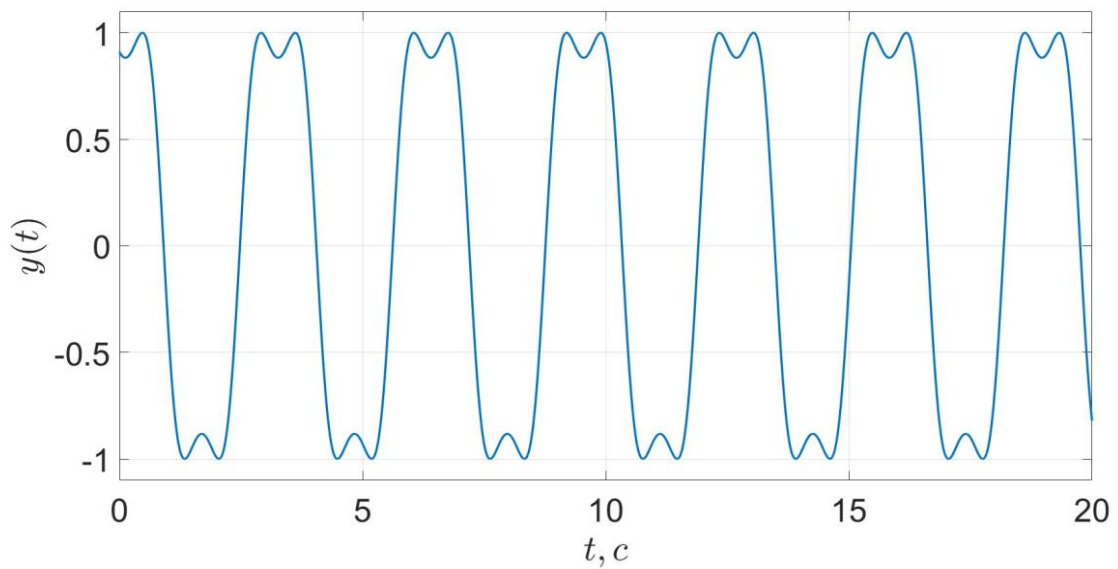

Рис. 7 – График сигнала $y(t)$ при $\theta_1 = 2, \theta_2 = 1, \gamma = 4$.

## 5. Заключение

В работе предложен новый подход к решению задачи оценивания нестационарной частоты синусоидального сигнала (1) при условии, что амплитуда также является переменной функцией времени. Данная задача была решена при выполнении допущений вида (2), (3) и (4), (5).

В статье представлен пример, разъясняющий на конкретном случае процедуру синтеза алгоритма оценивания, а также приведены результаты компьютерного моделирования, иллюстрирующие достижение заданной цели для различных параметров модели изменения частоты.

В качестве дальнейшего развития предложенного результата видится его расширение на случай неизвестных матриц уравнений (3) и (5).

СПИСОК ЛИТЕРАТУРЫ

А.А. Бобцов, *Университет ИТМО, профессор, Санкт-Петербург,* bobtsov@mail.ru

Н.А. Николаев, *Университет ИТМО, доцент, Санкт-Петербург,*
nikona@yandex.ru

О.В. Оськина, *Университет ИТМО, аспирант, Санкт-Петербург,*
olga.oskina1996@gmail.com

С. И. Низовцев, *Университет ИТМО, аспирант, Санкт-Петербург,*
nizovtsev.si@gmail.com